\documentclass[11pt,palatino]{article}
\usepackage[margin=0.75in]{geometry}

\usepackage[numbers,sort&compress,square]{natbib}

\usepackage{amsthm, amsmath, mathtools, mathrsfs}

\usepackage{amsfonts,amsmath,amssymb,bm,url}

\usepackage{graphicx,epsfig,subfigure}

\usepackage[utf8]{inputenc} 
\usepackage[T1]{fontenc}    
\usepackage{hyperref}       
\usepackage{url}            
\usepackage{graphicx}
\usepackage{booktabs}       
\usepackage{amsfonts}       
\usepackage{nicefrac}       
\usepackage{microtype}      
\usepackage{xcolor}         
\usepackage{comment}
\usepackage{amsmath}
\usepackage{amsthm}
\usepackage[utf8]{inputenc}
\usepackage[english]{babel}
\usepackage{authblk}
\usepackage{algorithm}
\usepackage{algpseudocode}

\begin{document}

\title{Koopman Learning with Episodic Memory}

\author[1*]{William T. Redman}
\author[1]{Dean Huang}
\author[1]{Maria Fonoberova}
\author[1$\dagger$]{Igor Mezić}
\affil[1]{\small AIMdyn Inc.}
\affil[*]{redmanw@aimdyn.com}
\affil[$\dagger$]{mezici@aimdyn.com}
\date{}
\maketitle

\begin{abstract}%
Koopman operator theory has found significant success in learning models of complex, real-world dynamical systems, enabling prediction and control. The greater interpretability and lower computational costs of these models, compared to traditional machine learning methodologies, make Koopman learning an especially appealing approach. Despite this, little work has been performed on endowing Koopman learning with the ability to leverage its own failures. To address this, we equip Koopman methods---developed for predicting non-autonomous time-series---with an episodic memory mechanism, enabling global recall of (or attention to) periods in time where similar dynamics previously occurred. We find that a basic implementation of Koopman learning with episodic memory leads to significant improvements in prediction on synthetic and real-world data. Our framework has considerable potential for expansion, allowing for future advances, and opens exciting new directions for Koopman learning.
\end{abstract}

\section{Introduction}
\label{section: Introduction}
Over the past two decades, Koopman learning \citep{mezic2021koopman} has emerged as a leading framework for modeling complex, real-world dynamical systems. The development of numerical methods to efficiently approximate spectral objects associated with the Koopman operator from data \citep{rowley2009dmd, williams2015edmd, arbabi2017ergodic, lusch2018deep, mezic2022numerical} has led to high-quality forecasting of challenging time-series, such as  Earth's magnetic field \citep{brunton2017chaos}, highway traffic patterns \citep{avila2020data}, climate indicators \citep{hogg2020exponentially}, and COVID-19 spreading \citep{mezic2023koopman}. The generalization of Koopman operator theory to include controlled dynamical systems \citep{proctor2016dynamic, proctor2018generalizing, korda2018linear, mauroy2020koopman} has led to state-of-the-art control of soft robotics \citep{bruder2019modeling, bruder2020data, bruder2021koopman, haggerty2023control}. The advancement of theoretical results \citep{mezic2005spectral, mauroy2013iso, lan2013linearization, arbabi2017ergodic, mezic2020spectrum}, linking Koopman spectral objects to state-space geometry---as well as the inherent linearity of the Koopman operator---has enabled the interpretation of learned models to an unprecedented degree. 

Unlike many traditional machine learning (ML) methods, such as convolutional neural networks, that are trained to extract informative features of the data (i.e., observable functions), Koopman learning identifies dynamical relationships between observables. Thus, while requiring appropriate observables to be chosen by a knowledgeable human or ML model \citep{li2017extended, lusch2018deep, otto2019linearly, yeung2019learning}, Koopman learning can discover spatio-temporal interactions that are crucial to understanding the behavior of the system under study.

\begin{figure*}[t]
    \centering
    \includegraphics[width=0.925\textwidth]{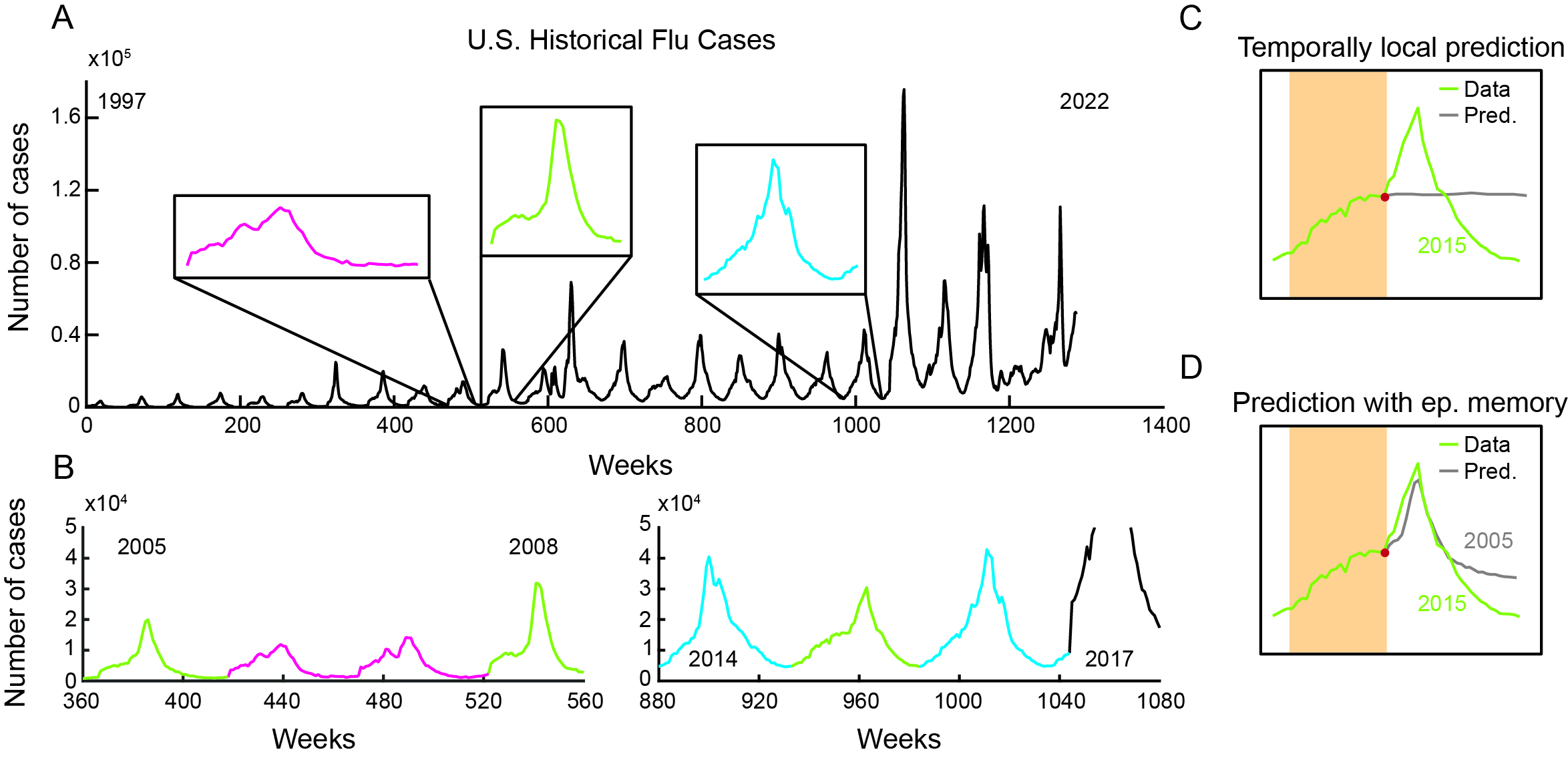}
     \caption{\small{\textbf{Example of repeated structure in non-autonomous dynamical system.} (A) Number of flu cases per week in the U.S., from 1997 to 2022, as reported to CDC by the Public Health Laboratories. (B) Zoomed-in four year time windows of (A). Years are colored by their qualitative similarities in shape and peak amplitude. (C) Predicting the future number of cases at a given time point (denoted by the red circle), using a temporally local window (denoted by an orange box), can lead to poor predictions (compare green and gray lines). (D) However, if a similar dynamical regime has occurred in the past, ``remembering'' this and using it for prediction can lead to improved performance (compare green and gray lines). }}
    \label{fig:Motivation schematic}
\end{figure*}

Koopman learning's data-driven, predictive, efficient, controllable, interpretable, and dynamics-based properties make it well aligned with next-generation artificial intelligence. And yet, to-date, there has been almost no work on enabling Koopman learning to leverage its failures, an essential aspect of any ``intelligent'' framework. The need for such an approach becomes clear when considering applying Koopman learning to non-autonomous dynamical systems [such as the yearly flu cases in the United States (U.S.), Fig.~\ref{fig:Motivation schematic}A], where temporally-local Koopman models are generated over streaming or sliding time windows \citep{macesic2018koopman, avila2020data, mezic2023koopman}. Although this mitigates the mixing of data from distinct dynamical regimes, the implicit forgetting of the past requires any repeated structure of the data (Fig.~\ref{fig:Motivation schematic}B) to be re-learned from scratch. 

Traditional ML methods have long wrestled with this problem, and---in their application to time-series---the use of memory has proven to be useful. This was done first by recurrent neural networks (RNNs) \citep{elman1990finding}, whose within-layer connectivity enables information to remain in the network through time, affecting future internal states. The introduction of long short-term memory (LSTM) networks \citep{hochreiter1997long} and liquid state machines \citep{maass2002real} led to computationally-tractable models that became widely adopted. However, in these models, long histories (or ``memory traces'') are continually present, making training challenging and susceptible to vanishing gradients. Over the past several years, Transformer networks \citep{vaswani2017attention}, which make use of the self-attention mechanism---originally designed to improve upon RNNs used for machine translation \citep{bahdanau2014neural}---have performed dominantly in a number of domains, attesting to the power \textit{selectively} attending to global aspects of data can provide. In the context of time-series prediction, where Transformers have become increasingly utilized \citep{wu2021autoformer, liu2022non, nie2022time, zhang2022crossformer}, self-attention can be viewed as recalling specific points in the past that are relevant to the present. Despite these advances, the high computational costs required to train such models, as well as the mixed performance gains reported when comparing Transformers to linear models \citep{zeng2023transformers}, suggests that there may be room for improvement from alternative approaches. 

Inspired by the integration of memory into traditional ML methods, we develop an approach that saves spectral objects associated with temporally local approximations of the Koopman operator, and utilizes this information to make new predictions that are directly informed from the past. Because the Koopman spectral objects are computed over a time window, these representations of the dynamics capture \textit{episodes} of previously observed data. Therefore, we refer to the general framework we introduce as \textbf{Koopman learning with episodic memory}. Through deployment on synthetic and real-world data sets, we show that the specific realization of Koopman learning with episodic memory we develop, called \textbf{extended dynamic mode decomposition (EDMD) with memory}, sees large improvements in prediction accuracy. This comes in large part from the ability of these models to not repeatedly make the same mistakes (Fig.~1C) by leveraging repeated structure in the data (Fig.~1D). We believe that our work not only demonstrates a way in which to better applications of Koopman operator theory to the prediction of non-autonomous dynamical systems, but additionally advances the state of Koopman learning as an artificial-intelligence framework, opening up a wealth of future directions. 

\section{Koopman operator theory}
The central object of interest in Koopman operator theory is the Koopman operator, $U$, an infinite-dimensional linear operator that describes the time evolution of observables (i.e., functions of the underlying state-space variables) that live in the functional space $\mathcal{F}$ (e.g. $\mathcal{L}^2$) \citep{Koopman1931, Koopman1932, mezic2005spectral}. In particular, after $t > 0$ amount of time, the value of the observable $f \in \mathcal{F}$ is given by
\begin{equation}
    \label{eq: KO}
    U^t f(x_0) = f \left[ T^t(x_0) \right],
\end{equation}
where $T$ is the dynamical map evolving the system and $x_0 \in \mathbb{C}^n$ is the initial condition in state-space. For the remainder of the paper, it will be assumed that the state-space being considered is of finite dimension and that $\mathcal{F}$ is a suitably chosen space of functions in which the spectral expansion exists \citep{mezic2020spectrum}.

In the case of a discrete-time dynamical system, the action of the Koopman operator on the observable $f$ can be decomposed as 
\begin{equation}
\label{eq: Koopman mode decomposition}
     U f(p) = \sum_{i = 0}^\infty \lambda_i \phi_i(x_0) v_i,
\end{equation}
where the $\phi_i$ are eigenfunctions of $U$, with $\lambda_i \in \mathbb{C}$ as their eigenvalues and $v_i$ as their Koopman modes \cite{mezic2005spectral, mezic2020spectrum}. For systems with chaotic or shear dynamics, Eq.~\ref{eq: Koopman mode decomposition} has an additional term corresponding to the continuous part of the spectrum \cite{mezic2005spectral}. For the remainder of the paper, we consider only the point spectra (Eq.~\ref{eq: Koopman mode decomposition}) of the associated Koopman operator.  

Spectrally decomposing the action of the Koopman operator is powerful because, for a discrete-time dynamical system, the value of $f$ at time step $k \in \mathbb{N}$ is given by
\begin{equation}
    \label{eq: Koopman mode decomposition time}
     f\left[ T^k(x_0) \right]  = U^k f(x_0) = \sum_{i = 0}^\infty \lambda_i^k \phi_i(x_0) v_i.
\end{equation}
\noindent
From Eq.~\ref{eq: Koopman mode decomposition time} it becomes clear that the dynamics of the system in the directions $v_i$, scaled by $\phi_i(x_0)$, are given by the magnitude of the corresponding $\lambda_i$. When presenting our algorithm below, we will (for simplicity) consider the output of the Koopman mode decomposition (KMD) being the Koopman eigenvalues, $\lambda_i$, and the scaled Koopman modes, $\phi_i(x_0)v_i$, which we will denote by $\tilde{v}_i(x_0)$. 

\subsection{Extended dynamic mode decomposition}
There exist many ways of numerically approximating the KMD from a data matrix $X = [x_0, x_1, ..., x_{m-1}]$, where $x_i \in \mathbb{C}^n$ are snapshots of the underlying dynamical system. In the case that a single observable (i.e., $n = 1$) is recorded from the system under study, $n_\text{delays} \in \mathbb{N}$, time delays can be introduced, so as to enlarge the set of observables \citep{takens2006detecting}. A popular approach to compute the KMD is dynamic mode decomposition (DMD) \citep{schmid2010dmd, rowley2009dmd}. DMD finds a matrix $K \in \mathbb{C}^{n \times n}$, such that 
\begin{equation}
\label{eq: DMD}
    || K X - Y ||_F,
\end{equation}
is minimized with respect to the Frobenius norm, where $Y = [x_1, ..., x_m]$. A solution to Eq.~\ref{eq: DMD} is given by 
\begin{equation}
    K = Y X^+,
\end{equation}
where $+$ is the Moore-Penrose pseudo-inverse \citep{rowley2009dmd, schmid2010dmd, tu2014dmd}. Although simple, this approach generates the KMD when certain conditions are met \citep{rowley2009dmd}. These conditions are restrictive, and in general, are not guaranteed to be satisfied. In addition, it is often the case that the action of the Koopman operator, restricted to the identity observable (that is, the original state-space variables), is not invariant. To address these limitations, extensions of DMD have been proposed.

One such approach is extended dynamic mode decomposition (EDMD) \cite{williams2015edmd}, which considers not the data matrices $X$ and $Y$, but a dictionary of functions $\mathcal{D} = [\psi_1, ..., \psi_l]$ applied to $X$ and $Y$. That is, it considers the functions $\Psi(x) = [\psi_1(x), ..., \psi_l(x)]$, where $\Psi:\mathbb{C}^n \rightarrow \mathbb{C}^l$. The Koopman operator is then approximated as
\begin{equation}
\label{eq: EDMD}
    K = G^+ A,
\end{equation}
where 
\begin{equation}
\centering
\begin{split}
        G = \frac{1}{m} \sum_{i = 1}^m [\Psi(x_i)]^*\Psi(x_i) \\
        A = \frac{1}{m} \sum_{i = 1}^m [\Psi(x_i)]^*\Psi(y_i).
\end{split}
\end{equation}

EDMD introduces a choice of which set of functions to include in the dictionary $\mathcal{D}$, which can be addressed from an understanding of the nature of the data \citep{bollt2021geometric} and/or can be optimized via dictionary learning \citep{li2017extended}. In this paper, we consider a family of Gaussian radial basis functions as the dictionary, but our approach generalizes beyond this decision. 

\subsection{Koopman mode decomposition for non-autonomous dynamical systems}

For non-autonomous dynamical systems, a Koopman mode decomposition can be rigorously defined \cite{mezic2016koopman, macesic2018koopman, proctor2018generalizing}, with the $N$ eigenvalues and scaled modes becoming time-dependent
\begin{equation}
\label{eq: non-autonomous Koopman mode decomposition}
    U g(x) = \sum_{i = 1}^N \lambda_i(t) \tilde{v}_i(x, t).
\end{equation}
A number of numerical methods for computing KMDs associated with non-autonomous dynamical systems have been developed \cite{kutz2016multiresolution, brunton2017chaos, arbabi2017ergodic, das2019delay, giannakis2019data, zhang2019online, avila2020data, ferre2023non, liu2023koopa, lu2023learning}. 

\begin{figure}[t]
    \centering
    \includegraphics[width=0.625\textwidth]{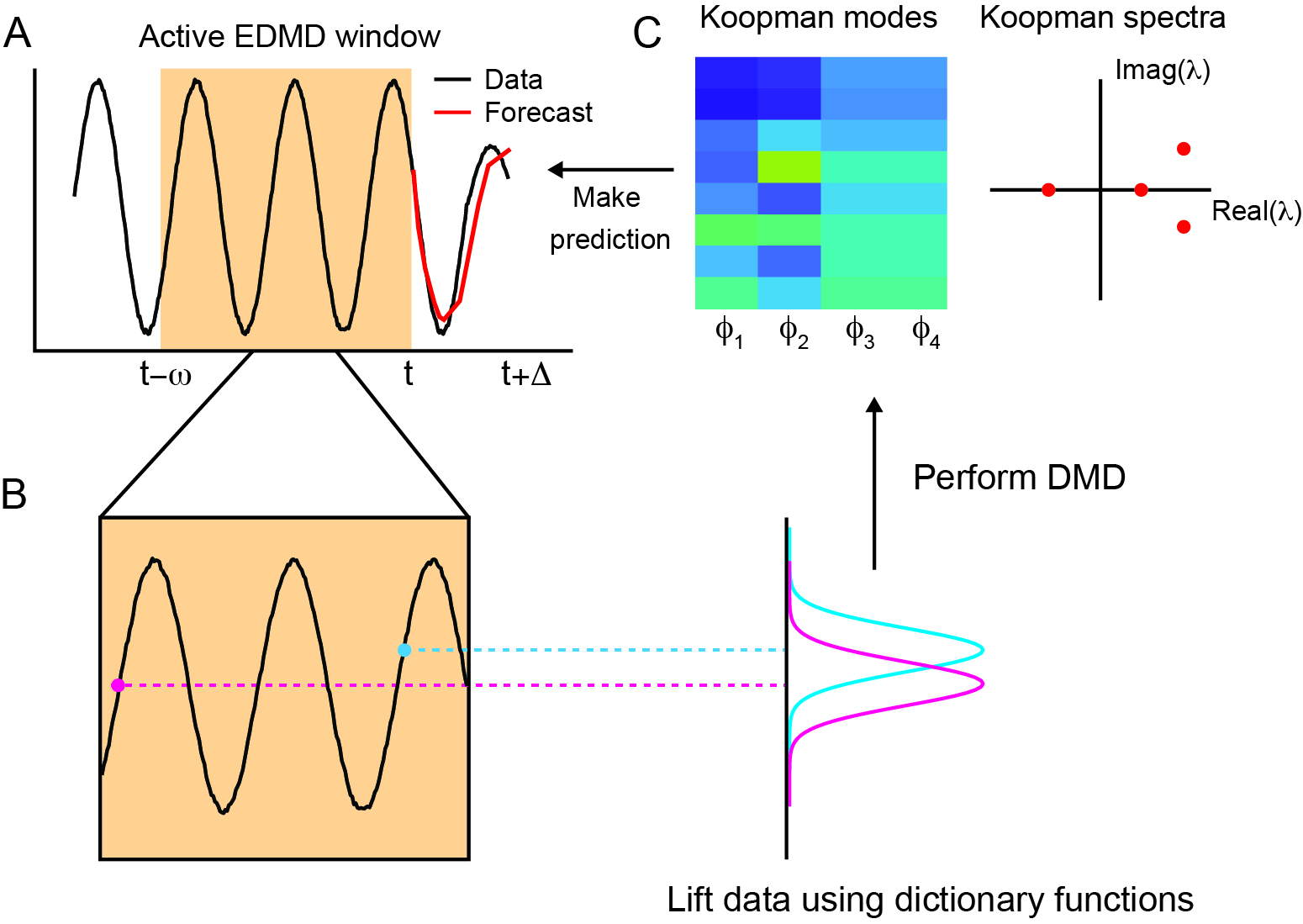}
     \caption{\small{\textbf{Illustration of sliding EDMD.}  Data (black line in A) from a temporally local window (orange box) is lifted (B), via a dictionary of functions (Eq.~\ref{eq: EDMD}---e.g., Gaussian radial basis functions). This lifting enables the approximation of a temporally local Koopman mode decomposition (C), which can then be used to make predictions $\Delta$ time steps into the future (red line in A). }}
    \label{fig: Sliding EDMD schematic}
\end{figure}

One simple approach is to learn temporally local Koopman operators over which the data can be approximated as autonomous. Let $X_{t, \omega} = [x_{t - \omega - 1}, ..., x_{t - 1}]$ and $Y_{t, \omega} = [x_{t - \omega}, ..., x_{t}]$ be $\omega \in \mathbb{N}$ snapshots ending at time step $t-1$ and $t$, respectively. A Koopman operator, $K_{t, \omega}$, can be computed using Eq.~\ref{eq: EDMD}, where $X_{t, \omega}$ and $Y_{t, \omega}$ are substituted for $X$ and $Y$. The Koopman mode decomposition (Eq.~\ref{eq: non-autonomous Koopman mode decomposition}) of $K_{t, \omega}$ yields Koopman eigenvalues, $\lambda_i(t, \omega)$, and the scaled Koopman modes, $\tilde{v}_i(t, \omega)$ (which, as a reminder, are the product of the modes scaled by the eigenfunctions). These can then be leveraged to predict $\Delta \in \mathbb{N}$ time steps into the future. This approach, referred to as ``sliding EDMD'' (Fig.~\ref{fig: Sliding EDMD schematic}), has been used for prediction in a number of applied settings \cite{brunton2017chaos, avila2020data, mezic2023koopman}, finding considerable success. Indeed, this algorithm led to state-of-the-art prediction of 1-week ahead COVID-19 forecasts \footnote{https://aimdyn.com/2023/06/13/forecasts-and-black-swans/}. This approach is similar to switched linear dynamical systems (SLDS) \cite{fox2008nonparametric, costa2019adaptive}, although SLDS's assumption of linearity in the state-space is not always met. In contrast, sliding EDMD can leverage a dictionary of nonlinear functions to approximate linear dynamics in an invariant subspace. 

\subsection{Comparing dynamical systems}
Given two dynamical maps, $T$ and $S$, it may be of interest to understand how ``similar'' their induced dynamics are. In the case where the dynamics of $T$ can be smoothly and invertibly mapped to the dynamics of $S$, the two are said to be topologically conjugate \cite{wiggins1996introduction}. That is, a topologically conjugacy exists if there exists a homeomorphism $h$, such that 
\begin{equation}
\label{eq: topological conjuagcy}
    h \circ T = S \circ h.
\end{equation}

Historically, identifying topological conjugacies, except in the simplest settings, has been challenging \cite{skufca2008concept}. However, recent theoretical work has proven that if $T$ and $S$ are in the basin of attraction of a stable fixed point, then they are topologically conjugate if and only if $\lambda^{(T)}_i = \lambda^{(S)}_i$ $\forall i = 1,...,N$, where $\lambda^{(T)}$ and $\lambda^{(S)}$ are the eigenvalues corresponding to the Koopman operators associated with $T$ and $S$ \cite{mezic2020spectrum}. This prior work utilized an extension of the Hartman-Grobman theorem to the entirety of the basin of stable fixed points \cite{lan2013linearization}, and does not hold when the system under study has a continuous Koopman spectrum. In such cases, it has been hypothesized that extensions of Koopman operator theory can be used to identify conjugacies \cite{avila2023spectral}. 

To quantify the similarity of Koopman eigenvalues, previous work \cite{redman2022algorithmic, redman2023equivalent, ostrow2024beyond} has utilized the Wasserstein distance \cite{kantorovich1960mathematical, vaserstein1969markov}, a metric developed in the context of optimal transport that quantifies how much one distribution must be changed to match the other. Formally, two probability measures $\mu$ and $\nu$, defined on the metric space $(M, d)$, have a $p^{\text{th}}$ Wasserstein distance given by 
\begin{equation}
\label{eq: Wasserstein distance}
  W_{p}(\mu ,\nu ):=\left(\inf _{\gamma \in \Gamma (\mu ,\nu )}\int _{M\times M}d(x,y)^{p}\,\mathrm {d} \gamma (x,y)\right)^{1/p},
\end{equation}
where $\Gamma(\mu, \nu)$ is the collection of all measures on $M \times M$, with marginals $\mu$ and $\nu$, respectively, and the $p^{\text{th}}$ moments of $\mu$ and $\nu$ are assumed to be finite. In the case where only a small, finite number of Koopman eigenvalues are computed, the Wasserstein distance is the distribution of eigenvalues, viewed as a sum of the Dirac measure. This enables the efficient computation of the Wasserstein distance by using linear sum assignment, making it highly scalable. Such a computational regime, which is the setting of the experiments performed in Sec.~\ref{sec: Results}, can be achieved by performing dimensionality reduction or residual-based pruning of modes when considering systems with a larger number of observables \cite{drmac2018data}.

When looking for periods in time when non-autonomous dynamical systems exhibit similar future behavior, identifying a topological conjugacy may not be enough. In particular, we expect it to be necessary that the Koopman modes are also equivalent, as this ensures that the same observables are evolving similarly. Therefore, in addition to computing the distance between pairs of Koopman spectra, we also compute the distance between pairs of Koopman modes. A simple way to achieve this is via the Euclidean distance between the modes. To make this independent of the norm of the modes, we can consider the normalized Euclidean distance. 

\section{Koopman learning with episodic memory}
\begin{figure*}[t]
    \centering
    \includegraphics[width = 0.925\textwidth]{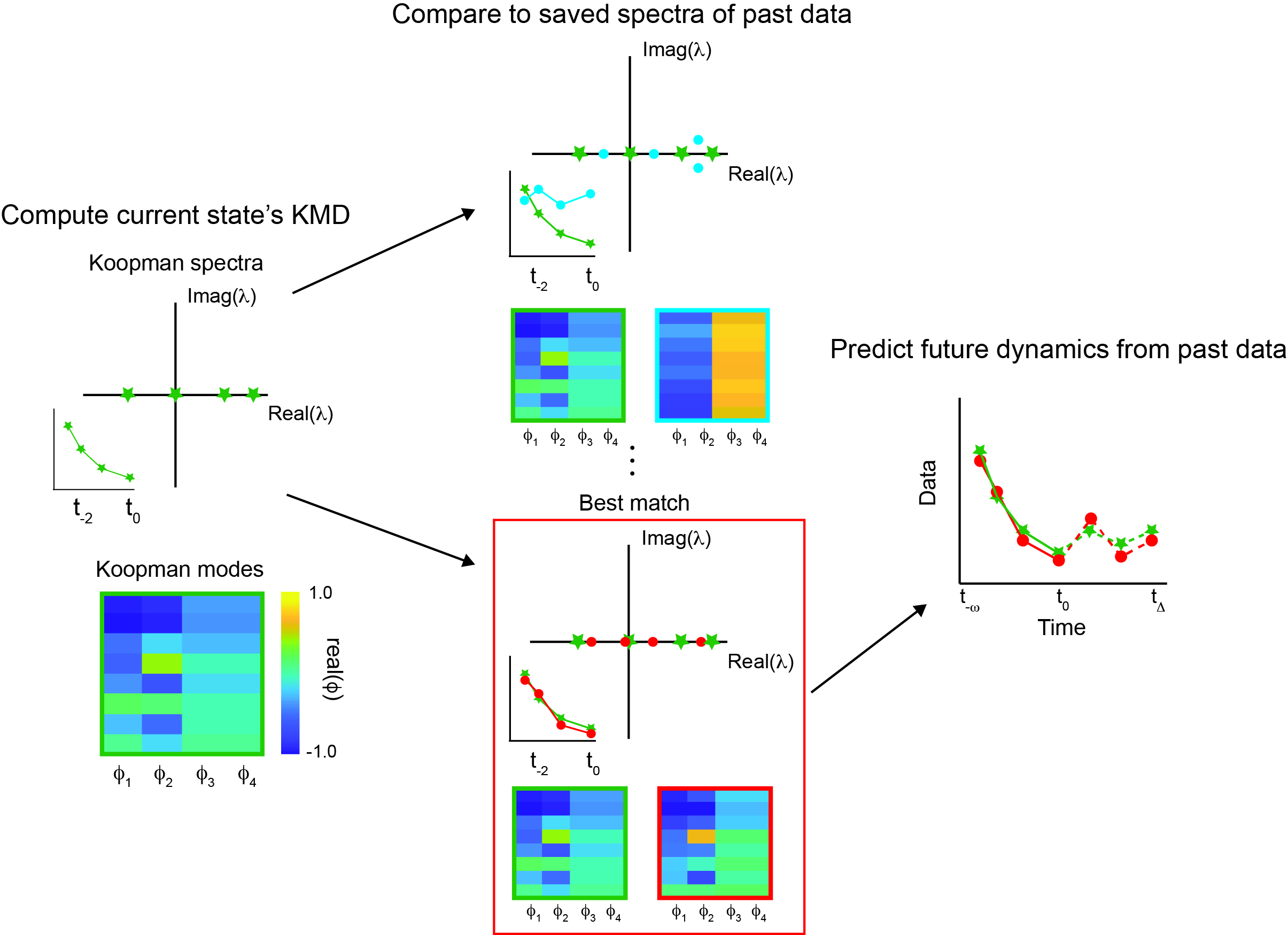}
    \caption{\small{\textbf{Overview of episodic memory for Koopman learning.} For the current state of the system under study, the data within the sliding window (inset figure on left) is used to compute the KMD. The resulting Koopman eigenvalues and modes are then compared to  a memory bank of saved Koopman eigenvalues and modes, associated with past data. If a match of both Koopman eigenvalues and modes is found (illustrated by the red box), then the data that occurred in the past (dashed red line) is used to predict the future behavior of the system (dashed green line).  If no match is found, sliding EDMD is used to generate a prediction.}}
    \label{fig:EDMD with memory schematic}
\end{figure*}

As noted in Sec.~\ref{section: Introduction}, a motivation for the integration of episodic memory in Koopman learning is the existence of repeated structure in some non-autonomous dynamical systems (Fig.~\ref{fig:Motivation schematic}B). These patterns can emerge when underlying interactions between components of the system repeat (e.g., human contact patterns), as well as when constraints of the system repeat (e.g., healthcare infrastructure and supply chains).

In order to leverage this repeated structure for improved prediction, one must identify when currently observed data has similar dynamical features to data seen in the past. This requires three ingredients. First, a representation of the dynamical features of data is needed. Given that Koopman operator theory decomposes data into Koopman eigenvalues and modes, the  representation of the KMD satisfies this requirement. In this paper, we use sliding EDMD to approximate the Koopman spectral objects, but our approach applies broadly to other numerical methods \citep{mezic2022numerical}. Second, a way to compare dynamical features is needed, so that the notion of ``similarity'' can be quantified. Although Koopman-based identification of topological conjugacy for non-autonomous dynamical systems is an underdeveloped area of research, we reason that it may be useful in this context. For this reason, we make use of the Wassertsein distance between the Koopman eigenvalues, and the Euclidean distance between the Koopman modes. And third, the dynamical features need to be saved for future comparison. Here, we simply store the Koopman eigenvalues and modes associated with \textit{all} local time windows in very basic data structures. Future applications of this approach can make use of more sophisticated ways in which to efficiently save the Koopman representations. 

\subsection{Algorithm}  

Our implementation of Koopman learning with episodic memory (\textbf{EDMD with memory} $-$ Fig.~\ref{fig:EDMD with memory schematic}, Algorithm \ref{algorithm:EDMD with memory}), on a data set $X = [x_0, ..., x_{m-1}]$ and $Y = [x_1, ..., x_{m}]$, has five primary steps:

\begin{enumerate}
    \item \textbf{Performing sliding EDMD:} At time step $t \geq \omega$, EDMD (Eq.~\ref{eq: EDMD}) is performed on $X_{t, \omega}$ and $Y_{t, \omega}$, using the dictionary $\mathcal{D} =  [\psi_1, ..., \psi_l]$. 

    \item \textbf{Computing distance to saved dynamics:} The computed Koopman eigenvalues and scaled modes, $\lambda(t, \omega)$ and $\tilde{v}(t, \omega)$, are compared, via the Wasserstein distance and the Euclidean norm, to the Koopman eigenvalues and modes saved from all previous time windows, $\lambda(t', \omega)$ and $\tilde{v}(t', \omega)$ s.t. $t' +\Delta < t$. In particular, $d_{\lambda}(t') = W_1[\lambda(t, \omega), \lambda(t', \omega)]$ (Eq.~\ref{eq: Wasserstein distance}) and $d_{v}(t') = \frac{||\tilde{v}(t, \omega) - \tilde{v}(t', \omega)||_2}{||\tilde{v}(t, \omega)||_2 + ||\tilde{v}(t', \omega)||_2}$ are computed.

    \item \textbf{Identifying a ``match'':} Let $t_\text{min} = \min_{t'} \{d_{\lambda}(t') + d_{v}(t') \big| d_{\lambda}(t') < \varepsilon_\lambda, \hspace{1mm} d_{v}(t') < \varepsilon_v \}$, where $\varepsilon_\lambda, \varepsilon_v \in \mathbb{R}^+$ are pre-chosen values. That is, $t_\text{min}$ is the time step in the past whose Koopman mode decomposition is closest to the Koopman mode decomposition at the current time, while also satisfying the condition that the distance between Koopman eigenvalues and modes are less than $\varepsilon_v$ and $\varepsilon_v$, respectively. If such a $t_\text{min}$ exists, then a dynamical ``match'' is said to be found. We find the best results in applying EDMD with memory to real-world data (Sec.~\ref{sec: Results}) when setting $\varepsilon_\lambda$ to be not too small (otherwise, no matches are found) and $\varepsilon_v$ to be proportional to the number time-delays used in the lifting. Concretely, we find that $\varepsilon_\lambda = 0.05$ or $0.10$ and $\varepsilon_v = \varepsilon_\lambda \cdot (n_\text{delays} + 1)$ works well.  

    \item \textbf{Predicting (with episodic memory):}  If a match is found, then a prediction $\Delta$ steps into the future is made by directly using what previously occurred $\Delta$ time steps after $t_\text{min}$. That is, we set $\hat{y}_{t + \Delta} = y_{t_\text{min} + \Delta}$, where $\hat{y}_{t + \Delta}$ is the prediction and $y_{t_\text{min} + \Delta}$ is the value of the data that occurred after the previous time window with similar dynamics. In this way, EDMD with memory ``remembers'' what happened the last time similar dynamics occurred. If no match is found, then the standard sliding EDMD prediction is used. 

    \item \textbf{Storing Koopman eigenvalues and modes:} After making a prediction, $\lambda(t, \omega)$ and $\tilde{v}(t, \omega)$ are stored with all the other previously saved spectral objects. 
    
\end{enumerate}

\begin{algorithm}[ht]
\caption{EDMD with memory.}
\label{algorithm:EDMD with memory}
\small{\begin{algorithmic}
\State \textbf{Input:} $X \in \mathbb{C}^{n \times m}$, $Y \in \mathbb{C}^{n \times m}$, $\Delta \in \mathbb{N}$, $\omega \in \mathbb{N}$, $\varepsilon_\lambda \in \mathbb{R}^+$, $\varepsilon_v \in \mathbb{R}^+$, and $\mathcal{D}$.
\State \textbf{Initialize:} V = [\hspace{1mm}], $\Lambda = [\hspace{1mm}]$, $\hat{Y} = [\hspace{1mm}]$
\For{$t = \omega + 1,..., m - \Delta - 1$}
\State $X_{t, \omega} = [x_{t - \omega - 1}, ..., x_{t - 1}]$ \small{(Step 1)} 
\State $Y_{t, \omega} = [x_{t - \omega}, ..., x_{t}]$ \small{(Step 1)} 
\State $\lambda(t, \omega), \tilde{v}(t, \omega) = \texttt{EDMD}(X_{t, \omega}, Y_{t, \omega}, \mathcal{D})$ \small{(Step 2)}
\State $t_\text{min} = 0$
\State $d = 100$
\For{$t' = \omega + 1, ..., t - \Delta$}
\State $d_\lambda = W_1(\lambda(t, \omega), \Lambda_{t'})$ 
\State $d_v = \frac{||\tilde{v}(t, \omega) - V_{t'}||_2}{||\tilde{v}(t, \omega)||_2 + ||V_{t'}||_2}$ 
\If{$d_\lambda < \varepsilon_\lambda$ and $d_v < \varepsilon_v$}
\If{$d > (d_\lambda + d_v)$}
\State $d = d_\lambda + d_v$
\State $t_\text{min} = t'$ \small{(Step 3)}
\EndIf
\EndIf
\EndFor
\If{$t_\text{min} > 0$}
\State $\hat{y} = Y_{t_\text{min} + \Delta}$ \small{(Step 4)}
\State $\hat{Y} = [\hat{Y}, \hat{y}]$
\Else 
\State $\hat{y} = \sum_i \lambda_i^\Delta(t, \omega) \cdot \tilde{v}_i (t, \omega)$ \small{(Step 4)}
\State $\hat{Y} = [\hat{Y}, \hat{y}]$
\EndIf
\State $\Lambda = [\Lambda, \lambda(t, \omega)]$ \small{(Step 5)}
\State $V = [V, \tilde{v}(t, \omega)]$ \small{(Step 5)}
\EndFor
\end{algorithmic}}
\end{algorithm}

Note that in our current implementation, we only allow information from one memory to be used, and for the prediction to be exactly what occurred before (up to scaling, see Sec.~\ref{sec: Results}). Combining information across multiple similar memories is an exciting direction that we will pursue in future work. 

\subsection{Computational complexity}

EDMD with memory has two primary computational costs. First, the temporally-local Koopman operator must be approximated from the sliding window. Let $l = |\mathcal{D}|$ be the number of functions in the dictionary. Then from Eq.~\ref{eq: EDMD} the computational complexity for approximating $K_{t, \omega}$ is $\mathcal{O}(l^3)$. 

This first cost is shared with any sliding EDMD-like method, and is independent of the number of time-windows. The second cost is specific to our episodic memory implementation, and is the cost of computing the Wasserstein distance between $\lambda(t, \omega)$ and all other saved eigenvalues. Because we consider a small number of dictionary elements ($l$ is assumed to be small), we can use linear sum assignment, which has computational complexity $\mathcal{O}(l^3)$. The total cost grows linearly with the number of time-windows that have been seen thus far ($t - \Delta - 1$). When $t$ is large, the total computational cost is therefore $\mathcal{O}(tl^3)$. 

To make this approach more efficient, a buffer of fixed sized can be used for the saved Koopman eigenvalues and modes, with $\lambda(t', \omega)$ and $\tilde{v}(t', \omega)$ being removed if $t - t'$ is large and/or if $\lambda(t', \omega)$ and $\tilde{v}(t', \omega)$ have not been recently matched with observed data. However, despite the linear cost with number of observed windows, in practice we had no problem running our unoptimized implementation on a laptop with a single CPU in minutes to (at most) an hour. 

\section{Results}
\label{sec: Results}
\subsection{Synthetic piecewise exponential data}

\begin{figure*}[t]
    \centering
    \includegraphics[width = 0.85\textwidth]{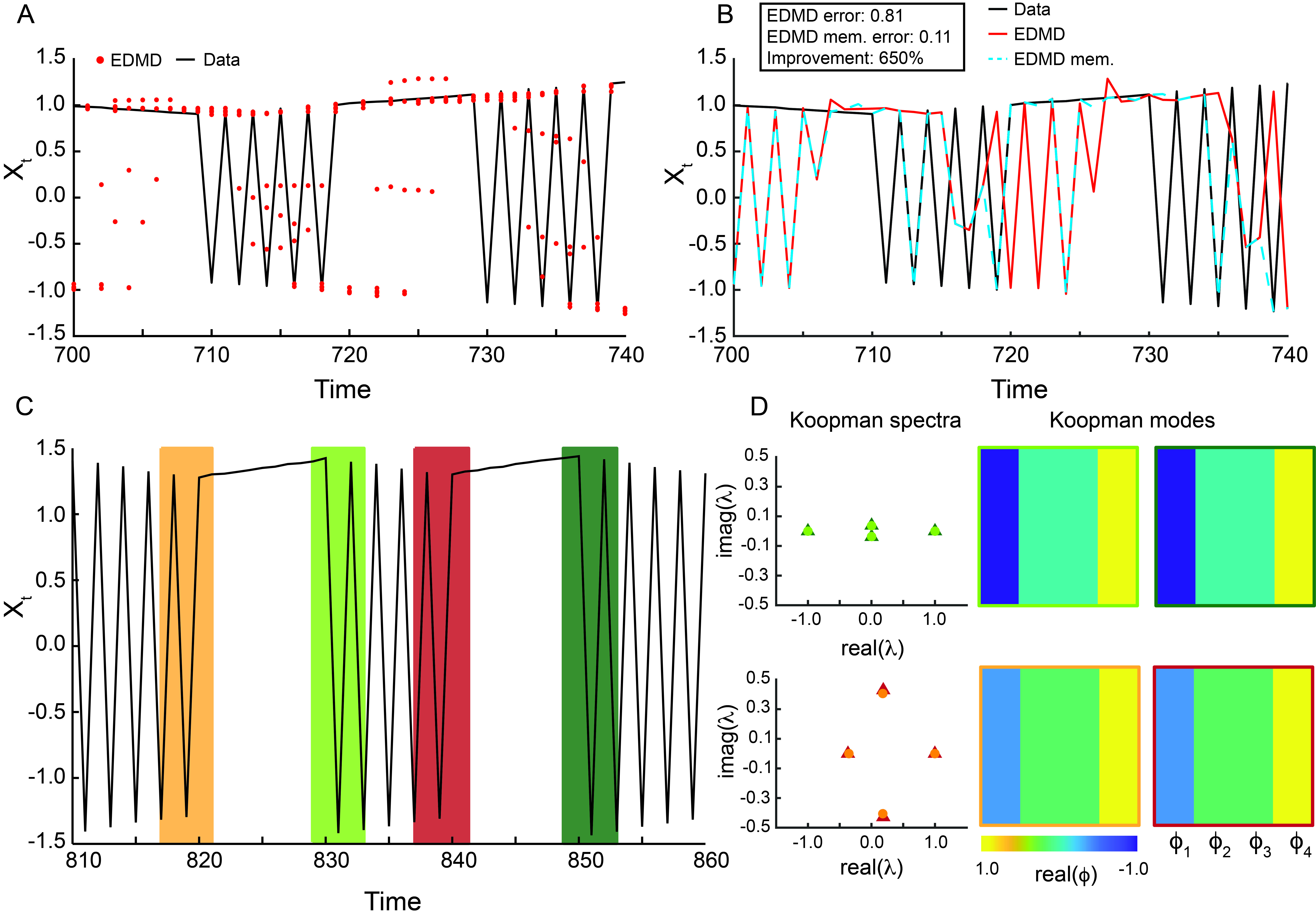}
    \caption{\small{\textbf{Episodic memory improves the ability of Koopman learning to predict synthetic non-autonomous exponential data.} (A) The predictions of sliding EDMD $\Delta = 1, ..., 5$ time steps ahead (red dots) often poorly match the true data (black line). (B) Utilizing EDMD memory (blue line) leads to $650\%$ improvement on predictions, as compared to sliding EDMD (red line). Note that although a zoomed-in section is shown here, the improvement was computed over a $1000$ time-step simulation. Here $\Delta = 5$, $\omega = 5$,  $n_\text{delays} = 1$, $\varepsilon_\lambda = 0.05$, and  $\varepsilon_v = 0.10$. (C) Because there are a finite (and small) number of possible values $\lambda_t$ can take, there is repetition of transitions (examples highlighted by light green/dark green and orange/dark red rectangles). (D) These repetitions can be identified by characteristic Koopman spectral objects [colored the same as in (C)], demonstrating how EDMD with memory is able to successfully recall relevant past data.}}
    \label{fig: non-autonomous noisy exponential}
\end{figure*}

To probe the predictive capabilities of EDMD with memory, we first designed a synthetic data set that standard sliding EDMD performs poorly on. In particular, we generated data coming from a piecewise exponential dynamical system described by 
\begin{equation}
    x_t = \lambda_t x_{t - 1} (1 + \eta \xi),
\end{equation}
where $\lambda_t \in \{-1.0101, -0.99, 0.99, 1.0101\}$ and is pseudo-randomly sampled every $10$ time steps\footnote{To mitigate numerical instability, we enforced the condition that the sign of $\lambda$ should change at each transition, so that the magnitude of $x_t$ does not explode.}, and $\xi$ is a random variable drawn from a Gaussian distribution of mean $0$ and variance $1$, with $\eta<<1$ acting as a magnitude of the noise.

As can be seen, the predictions from sliding EDMD do a poor job matching the true evolution of $x_t$ when the sliding window includes data from two different values of $\lambda_t$ (Fig.~\ref{fig: non-autonomous noisy exponential}A, B). However, although the Koopman models from these windows are not directly helpful for prediction, they do have identifiable Koopman eigenvalues and modes across repeats of the same transitions (e.g., $\lambda_t = -0.99 \rightarrow \lambda_t = 1.0101$) (Fig.~\ref{fig: non-autonomous noisy exponential}C, D). This suggests that utilizing the EDMD with memory may lead to improvements in prediction.

Because $x_t$ can grow and shrink with time (depending on whether $|\lambda_t| > 1$ or $|\lambda_t| < 1$), time windows with matching dynamics (and thus, similar Koopman eigenvalues and modes) can have different magnitudes of $x_t$. This can complicate directly using previously observed data as predictions. We correct for this by computing the ratio between the current time point's magnitude and that of the first time point in the recalled time window, and using this ratio as a re-scaling factor for our prediction. Additionally, we adaptively set the variance of the Gaussian radial basis dictionary functions to be the maximum value in the current sliding window, instead of the maximum value over the entire data set. This helps mitigate slight differences in Koopman spectra that artificially arise due to differences in scale. 

We find that the episodic memory mechanism can indeed boost the predictive capability of sliding EDMD on this data set (Fig.~\ref{fig: non-autonomous noisy exponential}B). While not perfect, EDMD with memory is more capable at accurately predicting the future state of $x_t$, even when the data in its sliding window comes from two different dynamical regimes (for instance, it can predict periods of oscillation, $\Delta = 5$ time steps ahead, from a sliding window that includes both exponentially growing and oscillating data). This suggests that the episodic memory mechanism is mitigating the repetition of the same mistakes, once a given transition has been observed. To quantify the improvement, we compute the median relative prediction improvement
\begin{equation}
  100\% \times \left[\frac{\text{median}_t(|x_{t + \Delta} - \hat{x}_{t + \Delta}^\text{EDMD}|)}{\text{median}_t(|x_{t + \Delta} - \hat{x}_{t + \Delta}^\text{mem}|)} - 1\right],
\end{equation}
where $\hat{x}_{t + \Delta}^\text{EDMD}$ is sliding EDMD's prediction $\Delta$ time steps ahead at time $t$, $\hat{x}_{t + \Delta}^\text{mem}$ is EDMD with memory's prediction $\Delta$ time steps at time $t$, and $\text{median}_t(\cdot)$ is the median value taken over all values of time. We find that EDMD with memory leads to an improvement of $650\%$ in median relative prediction error, as compared to sliding EDMD (Fig.~\ref{fig: non-autonomous noisy exponential}B). 

\subsection{U.S. Flu Cases}

We now turn our attention to applying EDMD with episodic memory on real-world, non-autonomous data. As mentioned in the Sec.~\ref{section: Introduction}, the number of U.S. flu cases per week\footnote{Data from https://gis.cdc.gov/grasp/fluview/fluportaldashboard.html.} has repeated structure (Fig.~\ref{fig:Motivation schematic}) that our episodic memory mechanism should be able to take advantage of. 

\begin{figure*}
    \centering
    \includegraphics[width = 0.999\textwidth]{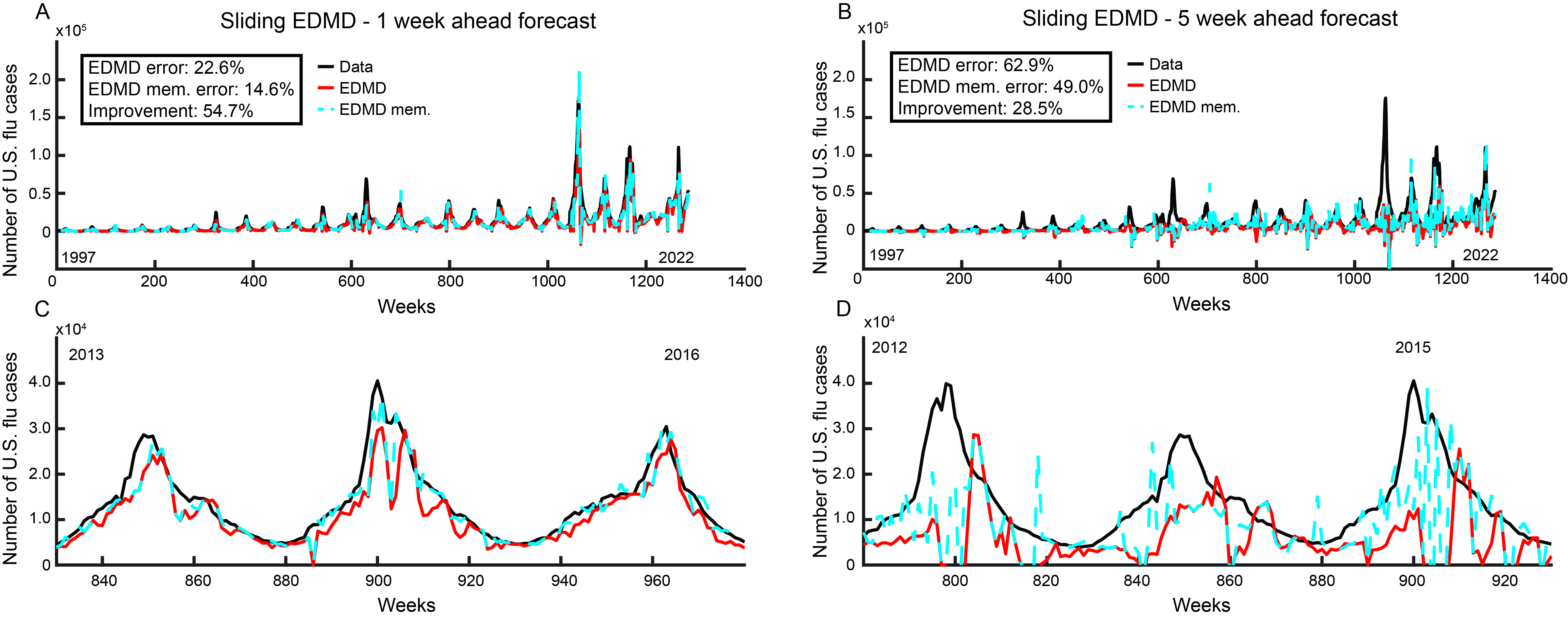}
    \caption{\small{\textbf{Episodic memory improves the ability of Koopman learning to predict U.S. flu cases.} (A)--(B) Predicted weekly number of U.S. flu cases one and five weeks ahead ($\Delta = 1$ and $5$, respectively), using sliding EDMD with (blue line) and without (red line) episodic memory. Black line denotes true data. (C)--(D) Zoomed-in windows of (A) and (B). For all subplots, $\omega = 3$,  $n_\text{delays} = 4$, $\varepsilon_\lambda = 0.05$, and  $\varepsilon_v = 0.25$.}}
    \label{fig:US flu cases}
\end{figure*}

We find that adding the episodic memory mechanism to sliding EDMD yields an improvement of almost $55\%$ on one-week ahead forecasting (Fig.~\ref{fig:US flu cases}A), decreasing the median relative error from $22.6\%$ to $14.6\%$. Examples of EDMD with memory's predictions across individual flu seasons (Fig.~\ref{fig:US flu cases}C) clearly show better matching with the ground truth. When increasing the forecast horizon to 5 weeks ahead, we again find that EDMD with memory can lead to large improvements over sliding EDMD (Fig.~\ref{fig:US flu cases}B, D), decreasing the error by over $28\%$. 

\begin{figure}
    \centering
    \includegraphics[width=0.875\linewidth]{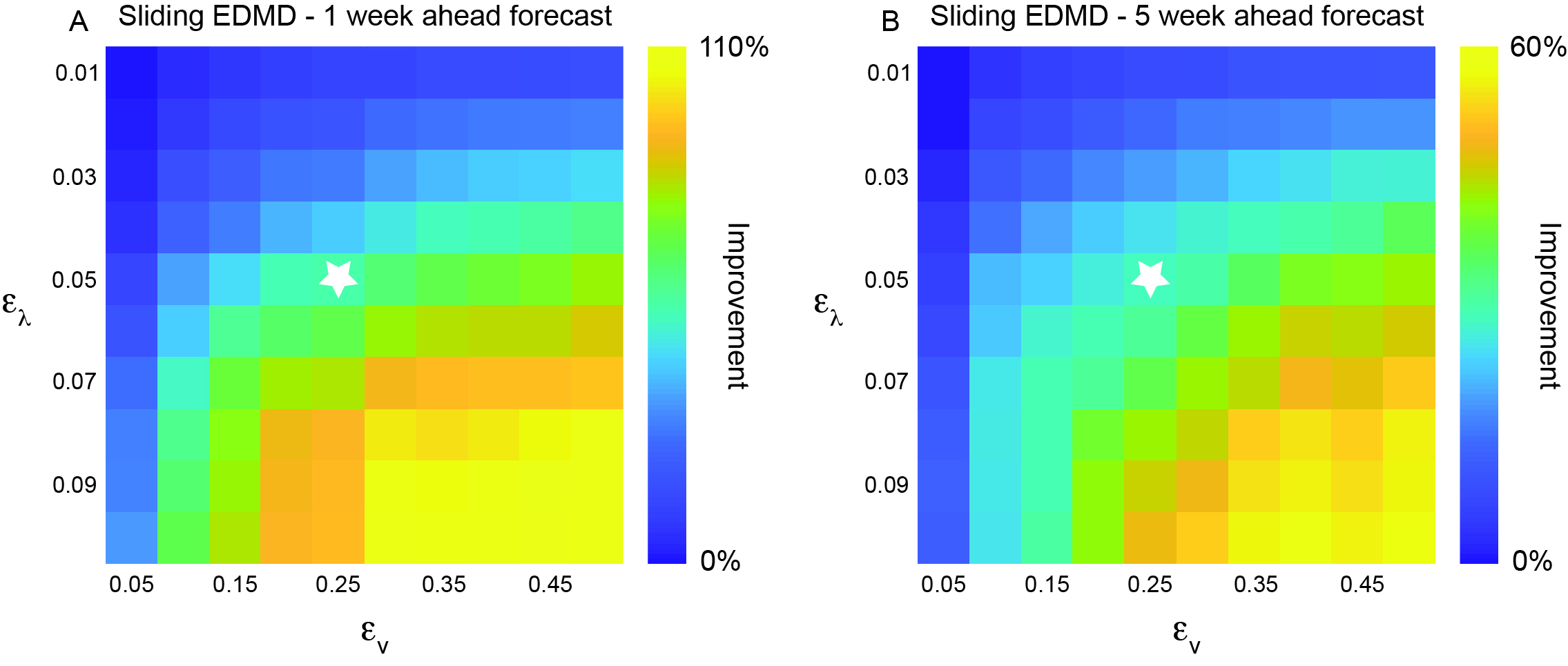}
    \caption{\small \textbf{EDMD with memory robustly provides benefit over sliding EDMD across a range of $\varepsilon_\lambda$ and $\varepsilon_v$}. (A) Improvement of using EDMD with memory for 1 week ahead forecasts of U.S. flu cases, as compared to sliding EDMD. (B) Same as (A), but for 5 week ahead forecasts. White star denotes values of $\varepsilon_\lambda = 0.05$ and $\varepsilon_v = 0.25$ used in Fig.~\ref{fig:US flu cases}. Note that for (A) and (B) the color bar does not have any negative values, demonstrating only observed improvement when using EDMD with memory. }
    \label{fig:US_flu_cases_epsilon_sweep}
\end{figure}

To ensure that these improvements are not due to ``lucky'' choices of $\varepsilon_\lambda$ and $\varepsilon_v$, we performed a grid search of $[\varepsilon_\lambda, \varepsilon_v] \in [0.01, 0.02, ...,0.1] \times [0.05, 0.10, ..., 0.50]$. We find a robust benefit of using EDMD with memory over standard sliding EDMD (Fig.~\ref{fig:US_flu_cases_epsilon_sweep}). As $\varepsilon_\lambda \rightarrow 0$ and $\varepsilon_v \rightarrow 0$, few memories are utilized, as the threshold for ``recalling'' becomes increasingly conservative. In this case, standard sliding EDMD is recovered and we find little improvement (Fig.~\ref{fig:US_flu_cases_epsilon_sweep}A$-$B, top left corners). In contrast, increasing $\varepsilon_\lambda$ and $\varepsilon_v$ allows for the incorporation of more memories and, for this data set, leads to better performance (Fig.~\ref{fig:US_flu_cases_epsilon_sweep}A$-$B, bottom right corners). The exact sensitivity to false positives (i.e., recalling a memory that is not relevant) depends on the application, and may require some tuning. However, for forecasting at both 1 and 5 weeks ahead on the U.S. flu cases, we find EDMD with memory only provides benefits, suggesting that it may be robust to the choice of $\varepsilon_\lambda$ and $\varepsilon_v$.

\subsection{Washington D.C. Bike Share}

\begin{figure*}
    \centering
    \includegraphics[width = 0.999\textwidth]{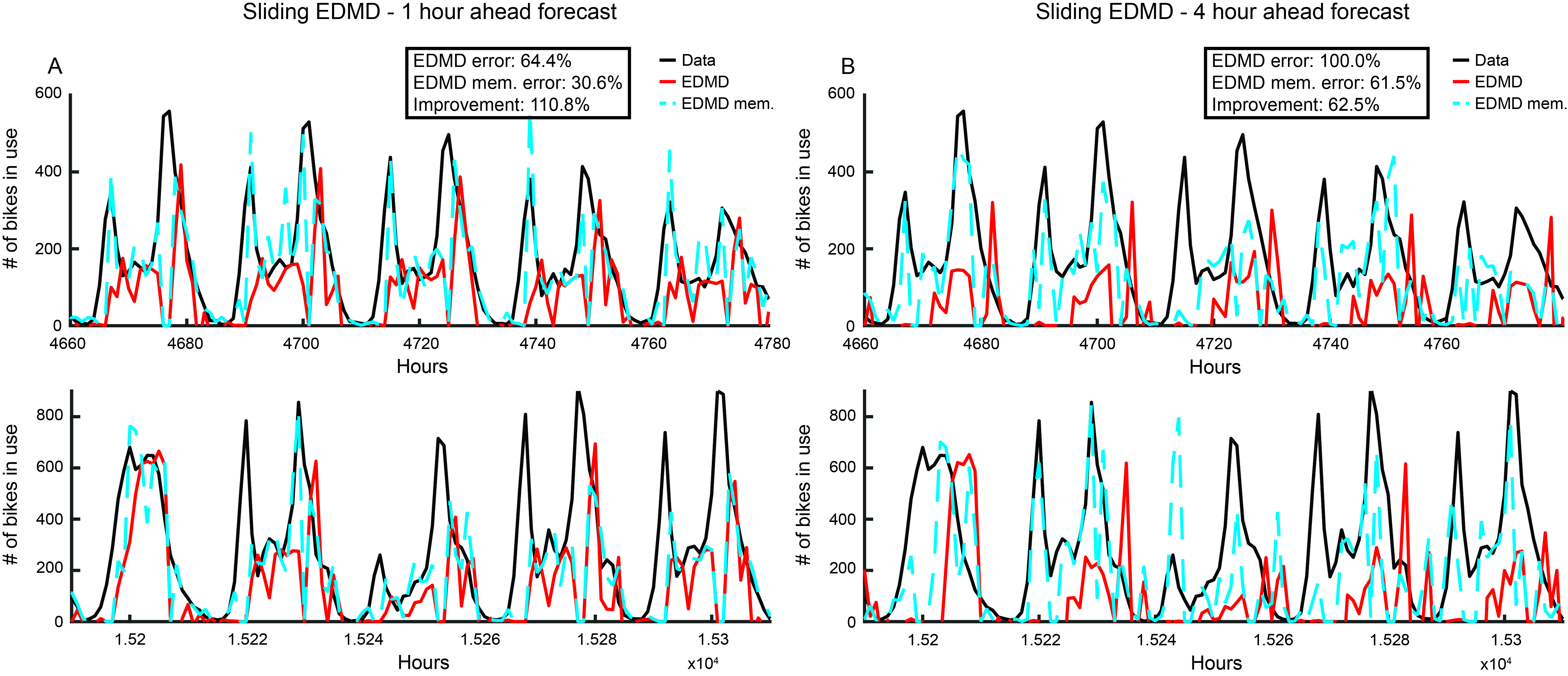}
    \caption{\small{\textbf{Episodic memory improves the ability of Koopman learning to predict bikes usage from the Washington D.C. bike share program.} (A) Predicted number of bikes used one hour ahead, from the Washington D.C. bike share program using sliding EDMD with (blue line) and without (red line) episodic memory. Due to the large number of time points ($17,000+$), the data has been zoomed-in around distinct five day (120 hour) periods (top and bottom). Note that the statistics presented in the top box are computed across the entirety of the data set. (B) Same as (A), but for four hour ahead predictions. For all subplots, $\omega = 3$,  $n_\text{delays} = 1$, $\varepsilon_\lambda = 0.1$, and  $\varepsilon_v = 0.2$.}}
    \label{fig:DC Bike Share Data}
\end{figure*}

Lastly, we apply EDMD with memory to predict the number of bikes in use, per hour, from the Washington D.C. bike share program during January 1, 2011 $-$ December 31, 2012 \footnote{Data from https://www.kaggle.com/datasets/marklvl/bike-sharing-dataset.}. This data set is especially challenging, given that it contains $17,000+$ time points and exhibits highly transient dynamics (Fig.~\ref{fig:DC Bike Share Data}A, B). As with the U.S. flu case data set, we expect there to be repeated dynamical structure due to repeated weather, holiday, and city traffic patterns. 

We find that EDMD with memory can provide improvements in median relative prediction error of more than $110\%$ over sliding EDMD, in the case of one-hour ahead forecasting, and over $62\%$, in the case of four-hour ahead forecasting (Fig.~\ref{fig:DC Bike Share Data}). We believe this improvement in performance, relative to the U.S. flu cases, comes from the bike share data having considerably more time points. This enables the episodic memory mechanism to leverage a wider repertoire of repeated patterns. Additionally, the variability in the data makes the application of sliding EDMD more challenging, allowing information ``recalled'' from the past, via the episodic memory mechanism, to have a larger impact.

\section{Discussion}
\label{section: Discussion}

In this paper, we present an episodic memory mechanism that can be easily integrated into the framework of Koopman learning (Fig.~\ref{fig:EDMD with memory schematic}) and can improve the predictive capabilities of sliding EDMD on non-autonomous data (Figs. \ref{fig: non-autonomous noisy exponential}--\ref{fig:DC Bike Share Data}). This is achieved by the ``recall'' of what occurred directly after episodes in the past that had similar dynamical features as those observed in the present. A consequence of this is that EDMD with memory is better able to avoid making the same mistakes each time it is given similar data. While our naive implementation of Koopman learning with episodic memory has limitations, it opens up a number of directions for expanding the ideas developed here that we believe will push forward the state of Koopman learning. 

\subsection{Related work}

As discussed in Sec.~\ref{section: Introduction}, this work was inspired in part by the use of attention in Transformers.  While our implementation of EDMD with memory is, like Transformers, able to selectively attend (or recall) previous episodes in time, there are two distinct differences between the two. First, unlike Transformers, EDMD with memory can be applied online, with no training. This greatly reduces computational costs, as illustrated by the fact that we were able to run EDMD with memory on the U.S. historical flu cases (Fig.~\ref{fig:US flu cases}) with a laptop using a single CPU and get results in minutes. Additionally, the online aspect of our approach enables continual integration of newly observed data, without requiring the re-training of a model. This greatly enhances its utility in applied settings. Second, we make use of KMD to represent dynamical information. This representation enables us to quickly search the past for episodes with similar dynamics, and make use of knowledge as to what happened. In contrast, Transformers must learn an embedding, which can be challenging. Indeed, initial attempts to apply Transformers to time series led to results that failed to outperform simple linear baselines \cite{zeng2023transformers}, with recent work suggesting that this may be due to Transformers learning sub-optimal embeddings \cite{ostrow2024delay}. Utilizing Koopman operator theory to inform the embeddings \cite{geneva2022transformers} has been shown to improve performance. This highlights the importance of the representation of the dynamics on which prediction algorithms are built, and it may be fruitful to combine EDMD with memory with Transformers utilizing Koopman embeddings to develop novel attention mechanisms. 

Our approach is also similar in spirit to approaches that attempt to extract KMDs for non-autonomous dynamical systems \cite{kutz2016multiresolution, ferre2023non, liu2023koopa, lu2023learning}. Such approaches face challenges not present in the autonomous case \cite{macesic2018koopman}, especially around times when dynamical regimes change. While EDMD with memory is unable to learn a complete representation of the non-autonomous dynamics, it is able to identify when previously observed data may be directly relevant, which can be useful at points where the dynamics change. Future work could seek to integrate the two types of approaches. 

Finally, the use of time delays and local linearization for short-term prediction of complex, and even chaotic, time-series has a long history \cite{eckmann1985ergodic, farmer1987predicting, casdagli1989nonlinear, jimenez1992forecasting, abarbanel1994predicting}. Prior work has developed methods for finding local associations between past time points and current data by identifying optimal linear mappings in a time delayed space. Additionally, time delays have been used for detection of anomalies in time-series by performing clustering in the time delayed space \cite{povinelli2003new}. Such approaches are similar to EDMD with memory as they leverage the known behavior of ``similar'' past data to make predictions about the future. However, by considering the Koopman mode decomposition on time delays (in addition to other observables), we believe our approach is better able to identify when past time points exhibit similar \textit{dynamics} to current data. Future work could compare these methods in greater detail to shed light on when one method may be expected to outperform another.

\subsection{Limitations}

Although the presence of repeating structure in time-series data (e.g., Fig.~\ref{fig:Motivation schematic}B) motivates the development of approaches that can identify and leverage them for future prediction, the success of our method is necessarily limited in several ways. 

First, the repeated structure that EDMD with memory identifies is related to constraints and response properties of the system under study. For instance, in the case of U.S. weekly flu cases, the repeated structure is related to healthcare infrastructure and supply chains, weather, and flu variants. When these are relatively consistent over the span of several years, EDMD with memory can be used to identify past data that is directly relevant. However, in cases where the constraints rapidly evolve with time, we expect EDMD with memory to be less useful. 

A second limitation has to do with the fact that EDMD with memory requires having seen a given dynamical feature in the past to predict it. Although it can, in principle, learn the pattern in a single shot (making it considerably more efficient than some standard ML methods), during the first exposure it does not contribute to the prediction. Because sliding EDMD has been found to be a powerful tool for data-driven prediction \cite{avila2020data, mezic2023koopman}, this should not be severely limiting to the predictive capability of our approach. However, there may be other methods, such as physics--informed (or, more generally, domain--informed) techniques, that would be better to deploy in the case where first identification is critical. However, we note that the absence of a good match between the current set of spectral objects and the library of saved spectral objects, provides information. Namely, that a new regime, unlike what has been seen previously, is being entered. If the library is created from a sufficiently large set of time points, this could provide confidence that novel dynamics are being encountered, and appropriate measures could be taken.

Third, in our current implementation of EDMD with memory, all saved Koopman eigenvalues and modes are given equal weighting, and only one, of possibly several close (with respect to the Wasserstein and Euclidean distance) previously saved memories, are used for prediction. Because this approach does not take into account whether this past data is nearby in time to the current data and/or whether this estimate came from data with high noise or uncertainty, we imagine that this approach is sub-optimal and can be improved upon.

\subsection{Expansions}

There are several ways in which this preliminary development of Koopman learning with episodic memory can be expanded upon. First, recent developments in Koopman operator theory for reduced order modeling allow for the noisy components of the Koopman spectra to be removed and the uncertainty in the estimate to be quantified \citep{mohr2022koopman}. Integrating this into EDMD with memory will enable better estimates of the underlying Koopman eigenvalues, as well as better matching between saved spectral objects. Second, a powerful capability of Koopman operator theory is its natural connection to observables beyond the state-space variables. As many state-space variables of interest (e.g., flu cases) are naturally functions of other variables (e.g., weather), these additional observables can be added to better approximate the dynamics, in ways that existing predictive methods may struggle with. We imagine including more observables into estimates of the Koopman operator will make EDMD with memory more accurate. Third, developments in numerical KMD approximations beyond EDMD have shown promise \cite{kostic2022learning, kostic2024sharp}. Future work can explore the effect of integrating our episodic memory approach with these methods. Fourth, we believe it will be fruitful to make connections with methods that utilize control, such as neural episodic control \cite{pritzel2017neural} and online convex optimization with memory and non-stochastic control \cite{zhao2023non}. Finally, we believe that research into human memory and decision making \citep{kahneman1982judgment} can be connected with our work and used to inspire further advances in Koopman learning. 

\subsection*{Acknowledgements} We thank the members of AIMdyn Inc. for useful discussion on our results and Thomas Redman for pointing out \cite{kahneman1982judgment}. This material is based upon work supported by the DARPA Small Business Innovation Research Program (SBIR) Program Office under Contract No. W31P4Q-21-C-0007 to AIMdyn, Inc. This material is also based upon work supported by the Air Force Office of Scientific Research under award number FA9550-22-1-0531. Any opinions, findings and conclusions or recommendations expressed in this material are those of the authors and do not necessarily reflect the views of the DARPA SBIR Program Office or AFRL/AFOSR. Distribution Statement A: Approved for Public Release, Distribution Unlimited.

\small
\bibliographystyle{unsrt}
\bibliography{main.bib}

\end{document}